\documentclass[10pt]{article}
\usepackage{amsfonts}
\usepackage{latexsym}
\usepackage{cite}
\usepackage{amsmath,amsfonts,latexsym,amssymb}
\usepackage[mathscr]{eucal}
\usepackage{cases,color}
\usepackage{amsthm}

\usepackage[bf,small]{caption2}
\usepackage{float}
\usepackage{graphicx}
\usepackage{amsmath}
\usepackage{amssymb}
\usepackage[all]{xy}

\newtheorem{theorem}{theorem}[section]

\newtheorem{prop}[theorem]{Proposition}

\newtheorem{thm}[theorem]{Theorem}

\begin{document}

\title{\vspace{-2cm}\textbf{Integral structure of the skein algebra of the 5-punctured sphere}}
\author{\Large Haimiao Chen}
\date{}
\maketitle

\begin{abstract}
  We give an explicit presentation for the Kauffman bracket skein algebra of the $5$-punctured sphere over any commutative unitary ring.

\end{abstract}

\section{Introduction}

Let $R$ be a commutative ring with identity and a fixed invertible element $q^{\frac{1}{2}}$.
Let $\Sigma=\Sigma_{g,h}$ denote the $h$-punctured orientable surface of genus $g$.
The {\it Kauffman bracket skein algebra} of $\Sigma$ over $R$, denoted by $\mathcal{S}(\Sigma;R)$, is defined as the $R$-module generated by isotopy classes of framed links (which may be empty) embedded in $\Sigma\times[0,1]$ modulo the following {\it skein relations}:
\begin{figure}[h]
  \centering
  \includegraphics[width=8.5cm]{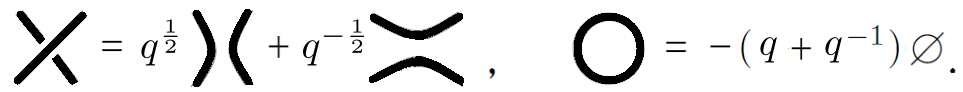}\\
\end{figure}

As conventions, $R$ is identified with $R{\O}\subset\mathcal{S}(M;R)$ via $a\mapsto a{\O}$, and each framed link is presented as a link with the understanding that each framing vector points vertically upward.

Given links $L_1,L_2$, the product $L_1L_2$ is defined by stacking $L_1$ over $L_2$ in the $[0,1]$ direction.
Using the skein relations, each element of $\mathcal{S}(\Sigma;R)$ can be written as a $R$-linear combination of multicurves on $\Sigma$, where a multicurve is identified with a link in $\Sigma\times\{\frac{1}{2}\}\subset\Sigma\times[0,1]$.

The description of the structure of $\mathcal{S}(\Sigma_{g,k};\mathbb{Z}[q^{\pm\frac{1}{2}}])$ was raised as \cite{Ki97} Problem 1.92 (J) and also \cite{Oh02} Problem 4.5. The structure of $\mathcal{S}(\Sigma_{g,k};\mathbb{Z}[q^{\pm\frac{1}{2}}])$ for $g=0,k\le 4$ and $g=1,k\le 2$ was known to Bullock and Przytycki \cite{BP00} early in 2000. A finite set of generators had been given by Bullock \cite{Bu99} in 1999. Till now it remains a difficult problem to find all relations for general $g$ and $k$.

As weak solutions, recently Cooke and Lacabanne \cite{CL22} obtained a presentation for $\mathcal{S}(\Sigma_{0,5};\mathbb{C}(q^{\frac{1}{4}}))$, and the author \cite{Ch24} found a presentation for $\mathcal{S}(\Sigma_{0,n+1};R)$ for all $n$ and all $R$ containing the inverse of $\alpha:=q+q^{-1}$.

In this paper we apply the result of \cite{Ch22} that the ideal of defining relations of $\mathcal{S}_n$ is generated by certain relations of degree at most $2n+2$, to the case $n=4$. We find an explicit set of relations to generate the ideal. The ``integral" in the title is justified by that $\alpha$ is no longer assumed to be invertible.
Actually our method can be easily extended to all $n\ge 5$.

Let $\overline{q}$ denote $q^{-1}$.
For a finite set $Y$, denote its cardinality by $\#Y$.

Recall some notations introduced in \cite{Ch22}; refer there for more details.

\begin{figure}[h]
  \centering
  \includegraphics[width=4cm]{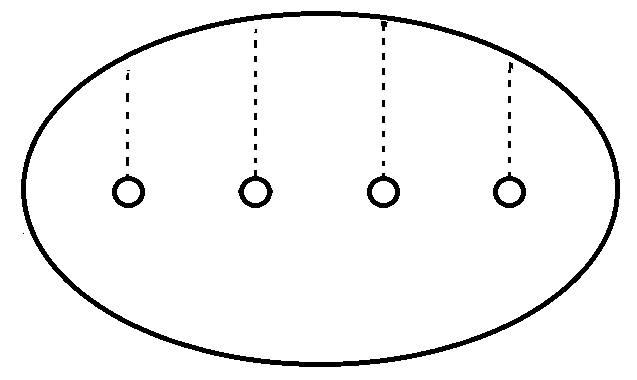}\\
  \caption{The surface $\Sigma_{0,4+1}$; the dotted lines are $\gamma_j$, $j=1,\ldots,4$.}\label{fig:Sigma}
\end{figure}

Display $\Sigma=\Sigma_{0,n+1}$ as in Figure \ref{fig:Sigma}. Let $\mathsf{p}_1,\ldots,\mathsf{p}_n$ denote the punctures, listed from left to right. Let $\gamma_i$ denote the vertical line connecting $\mathsf{p}_i$ and a point on $\partial\Sigma$. Let $\gamma=\bigcup_{i=1}^n\gamma_i$, and let $\Gamma=\bigcup_{i=1}^n\Gamma_i$, with $\Gamma_i=\gamma_i\times[0,1]$. Let $\pi:\Sigma\times[0,1]\to\Sigma$ be the projection.

Given a $1$-submanifold $X\subset\Sigma\times[0,1]$, define its {\it multidegree} as the function ${\rm md}_X:\{1,\ldots,n\}\to\mathbb{N}$, $v\mapsto\#(X\cap\Gamma_v)$.
Let $|X|=\#(X\cap\Gamma)$, called the {\it degree}.


For $1\le i_1<\cdots<i_r\le n$, let $t_{i_1\cdots i_r}\in\mathcal{S}_n$ denote the element represented by a simple curve exactly encircling $\mathsf{p}_{i_1},\ldots,\mathsf{p}_{i_r}$.
Let
$$\mathfrak{T}_n=\{t_{i_1\cdots i_r}\colon 1\le i_1<\cdots<i_r\le n,\ 1\le r\le n\}.$$
Denote $t_{1\cdots n}$ by $t_0$. Note that $t_1,\ldots,t_n$ and $t_0$ are central in $\mathcal{S}_n$.

Let $\mathcal{T}_n$ be the free $R$-algebra generated by $\mathfrak{T}_n$. A product of elements of $\mathfrak{T}_n$ is called a {\it monomial} and regarded as a link in $\Sigma\times[0,1]$.
Let $\theta_n:\mathcal{T}_n\to\mathcal{S}_n$ denote the canonical map.

For $\mathfrak{u}={\sum}_ia_i\mathfrak{g}_i\in\mathcal{T}_n$, where $0\ne a_i\in R$ and $\mathfrak{g}_i$ is a monomial, put
$${\rm md}_{\mathfrak{u}}(v)={\max}_i{\rm md}_{\mathfrak{g}_i}(v), 
\qquad |\mathfrak{u}|={\sum}_{v=1}^n{\rm md}_{\mathfrak{u}}(v).$$

Let $\Xi$ denote the set of $\vec{e}=(e_1,\ldots,e_n)$ with $\sum_{v=1}^ne_v\le 2(\#\{i\colon e_i>0\}+1)$. 
For each $\vec{e}\in\Xi$, let
$$\mathcal{R}(\vec{e})=\{\mathfrak{u}\in\ker\theta_n\colon{\rm md}_{\mathfrak{u}}(v)\le e_v,\ 1\le v\le n\}\subset\mathcal{T}_n.$$

Let $\mathcal{I}_n$ denote the two-sided ideal generated by the elements of $\bigcup_{\vec{e}\in\Xi}\mathcal{R}(\vec{e})$.

\begin{thm}[\cite{Ch22} Section 5]
The skein algebra $\mathcal{S}_{n}$ is generated by $\mathfrak{T}_n$, and the ideal of defining relations is $\mathcal{I}_n$.
\end{thm}

\section{A presentation for $\mathcal{S}_4$}

Suppose $J\subset\Sigma$ is a simple curve. 
Starting at a point $\mathsf{x}\in J$, walk along $J$ in any direction, record a label $i^\ast=i$ (resp. $i^\ast=\overline{i}$) whenever passing through $\gamma_i$ from left to right (resp. from right to left). Denote $J$ as $t_{i_1^\ast\cdots i_r^\ast}$ if when back to $\mathsf{x}$, the recorded labels are $i_1^\ast,\ldots,i_r^\ast$. This depends on the choices of $\mathsf{x}$ and the direction, so $J$ may have several different notations of such kind.

\subsection{Useful identities}

By direct computation,
\begin{align}
t_{12}t_{23}&=qt_{123\overline{2}}+\overline{q}t_{13}+t_1t_3+t_2t_{123}, \label{eq:12-23} \\
t_{13}t_{24}&=\alpha t_0+t_1t_{234}+t_2t_{134}+t_3t_{124}+t_4t_{123}+q^2t_{12}t_{34}+\overline{q}^2t_{14}t_{23} \nonumber \\
&\ \ \ \ +qt_3t_4t_{12}+\overline{q}t_1t_4t_{23}+qt_1t_2t_{34}+\overline{q}t_2t_3t_{14}+t_1t_2t_3t_4,   \label{eq:13-24}  \\
t_{14}t_{234}&=t_4t_0+\overline{q}t_{1\overline{4}234}+qt_{123}+t_1t_{23}, \label{eq:14-234} \\
t_{34}t_{124}&=t_4t_0+qt_{12\overline{4}34}+\overline{q}t_{123}+t_3t_{12}, \label{eq:34-124} \\
t_{24}t_{134}&=t_4t_0+\overline{q}t_{12\overline{4}34}+qt_{1\overline{4}234}+t_2t_{1\overline{4}34}.  \label{eq:24-134}
\end{align}

\begin{figure}[h]
  \centering
  \includegraphics[width=11cm]{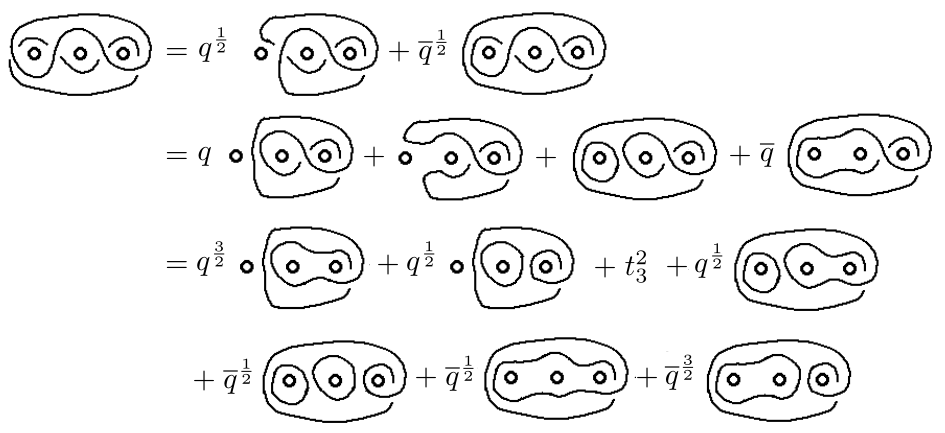}\\
  \caption{Computing $t_{123\overline{2}}t_{13}$.}\label{fig:top-times-down}
\end{figure}

Illuminated by Figure \ref{fig:top-times-down},
\begin{align}
t_{123\overline{2}}t_{13}&=t_{123}^2+(qt_1t_{23}+\overline{q}t_3t_{12}+t_1t_2t_3)t_{123}+q^2t_{23}^2+\overline{q}^2t_{12}^2 \nonumber \\
&\ \ \ \ +qt_2t_3t_{23}+\overline{q}t_1t_2t_{12}+t_1^2+t_2^2+t_3^2-\alpha^2.   \label{eq:top-times-down}
\end{align}

Since $t_{23}t_{34}=qt_{234\overline{3}}+\overline{q}t_{24}+t_2t_4+t_3t_{234}$, we have
\begin{align*}
t_{12}t_{23}t_{34}=q^2t_{1234\overline{3}\overline{2}}+t_{134\overline{3}}+qt_1t_4+qt_2t_{1234\overline{3}}+\overline{q}t_{12}t_{24}
+t_2t_4t_{12}+t_3t_{12}t_{234}.
\end{align*}
Then replacing $t_{134\overline{3}}$ and $t_{1234\overline{3}}$ by elements of $\mathcal{T}_n$, respectively through
\begin{align*}
t_{13}t_{34}&=qt_{134\overline{3}}+\overline{q}t_{14}+t_1t_4+t_3t_{134},   \\
t_{34}t_{123}&=\overline{q}t_{1234\overline{3}}+qt_{124}+t_4t_{12}+t_3t_0,
\end{align*}
we are led to
\begin{align}
t_{12}t_{23}t_{34}&=q^2t_{1234\overline{3}\overline{2}}+t_3t_{12}t_{234}+q^2t_2t_{34}t_{123}-q^3t_2t_{124}
-\overline{q}t_3t_{134}+\overline{q}t_{12}t_{24} \nonumber \\
&\ \ \ \ +\overline{q}t_{13}t_{34}+(1-q^2)t_2t_4t_{12}-\overline{q}^2t_{14}+(q-\overline{q})t_1t_4-q^2t_2t_3t_0. \label{eq:12-23-34}
\end{align}

\subsection{Relations}

\begin{prop}\label{prop:commuting}
In $\mathcal{S}_4$ we have
\begin{align}
qt_{23}t_{12}-\overline{q}t_{12}t_{23}&=(q^2-\overline{q}^2)t_{13}+(q-\overline{q})(t_1t_3+t_2t_{123}),  \label{eq:[2,2]-1} \\
t_{24}t_{13}-t_{13}t_{24}&=(\overline{q}^2-q^2)(t_{12}t_{34}-t_{14}t_{23})   \nonumber  \\
&\ \ \ +(\overline{q}-q)(t_3t_4t_{12}-t_1t_4t_{23}+t_1t_2t_{34}-t_2t_3t_{14}),  \label{eq:[2,2]-2} \\
\overline{q}t_{234}t_{14}-qt_{14}t_{234}&=(\overline{q}^2-q^2)t_{123}+(\overline{q}-q)(t_4t_0+t_1t_{23}), \label{eq:[2,3]-1}  \\
qt_{124}t_{34}-\overline{q}t_{34}t_{124}&=(q^2-\overline{q}^2)t_{123}+(q-\overline{q})(t_4t_0+t_3t_{12}), \label{eq:[2,3]-2} \\
t_{134}t_{24}-t_{24}t_{134}&=(q-\overline{q})(\overline{q}t_{34}t_{124}-qt_{14}t_{234}+qt_1t_{23}-\overline{q}t_3t_{12}) \nonumber  \\
&\ \ \ \ +(q-\overline{q})^2(t_4t_0+\alpha t_{123}).  \label{eq:[2,3]-3}
\end{align}
\end{prop}

\begin{proof}
The identities (\ref{eq:[2,2]-1})--(\ref{eq:[2,3]-2}) respectively follow from (\ref{eq:12-23})--(\ref{eq:34-124}) and their mirrors.
The last one follows from (\ref{eq:24-134}) and its mirror, as well as (\ref{eq:[2,3]-1}), (\ref{eq:[2,3]-2}).
\end{proof}

Call (\ref{eq:[2,2]-1}), (\ref{eq:[2,2]-2}) (resp. (\ref{eq:[2,3]-1})--(\ref{eq:[2,3]-3})) as well as the identities resulting from acting the indices via the permutations $(1234)^k,k=1,2,3$ the {\it commuting relations} of type $[2,2]$ (resp. type $[2,3]$).
They enable us to write any element of $\mathcal{S}_4$ as a linear combination of monomials in which the generators are arranged in any prescribed order.

Similarly as in \cite{Ch24}, we do not present commuting relations of type $[3,3]$.

\begin{prop}\label{prop:reduction}
The following relations hold in $\mathcal{S}_4$:
\begin{align}
t_{13}t_{24}&=\alpha t_0+t_1t_{234}+t_2t_{134}+t_3t_{124}+t_4t_{123}+q^2t_{12}t_{34}+\overline{q}^2t_{14}t_{23} \nonumber \\
&\ \ \ \ +qt_3t_4t_{12}+\overline{q}t_1t_4t_{23}+qt_1t_2t_{34}+\overline{q}t_2t_3t_{14}+t_1t_2t_3t_4,   \nonumber  \\
t_{24}t_{134}&=q^2t_{14}t_{234}+\overline{q}^2t_{34}t_{124}+(1-q^2-\overline{q}^2)t_4t_0-(q^3+\overline{q}^3)t_{123}-q^2t_1t_{23} \nonumber  \\
&\ \ \ \ -\overline{q}^2t_3t_{12}+t_2(qt_{14}t_{34}-q^2t_{13}-qt_4t_{134})-qt_1t_2t_3, \label{eq:24-134'}  \\
t_{123}^2&=\overline{q}t_{12}t_{23}t_{13}-(t_1t_2t_3+qt_1t_{23}+\overline{q}t_2t_{13}+\overline{q}t_3t_{12})t_{123}-(t_1^2+t_2^2+t_3^2) \nonumber  \\
&\ \ \ \ +\alpha^2-(qt_2t_3t_{23}+\overline{q}t_1t_3t_{13}+\overline{q}t_1t_2t_{12})-(q^2t_{23}^2+\overline{q}^2t_{13}^2+\overline{q}^2t_{12}^2), \label{eq:123-123}  \\
t_{123}t_{234}&=(t_{23}+qt_2t_3)t_0+\overline{q}t_{12}t_{23}t_{34}-\overline{q}t_3t_{12}t_{234}-qt_2t_{34}t_{234}+q^2t_2t_{124}  \nonumber  \\
&\ \ \ \ +\overline{q}^2t_3t_{134}-\overline{q}^2t_{12}t_{24}-\overline{q}^2t_{13}t_{34}+(q-\overline{q})t_2t_4t_{12}
+\overline{q}^2(\alpha t_{14}+t_1t_4),   \label{eq:123-234}  \\
t_{123}t_{134}&=t_{13}t_0+t_{12}t_{14}+t_{23}t_{34}-t_1t_{123}-t_3t_{234}-\alpha t_{24}-t_2t_4, \label{eq:123-134} \\
t_{23}t_{34}t_{124}&=(qt_{234}+q^2t_2t_{34}+t_3t_{24}+t_4t_{23}+qt_2t_3t_4)t_0 \nonumber \\
&\ \ \ \ +(t_2t_4+\overline{q}t_{24})t_{124}+q^3t_{34}t_{134}+\overline{q}t_{23}t_{123}+t_3t_{23}t_{12}+q^2t_{14}t_{34} \nonumber \\
&\ \ \ \ +qt_2t_{12}-q^3t_3t_{13}+qt_4t_{14}-q^2t_1t_3^2+q\alpha t_1, \label{eq:23-34-124} \\
t_{14}t_{12}t_{23}t_{34}&=q^2t_0^2+q^2(\overline{q}t_1t_{234}+qt_4t_{123}+t_1t_4t_{23}-t_2t_3)t_0+t_{234}^2+q^4t_{123}^2  \nonumber  \\
&\ \ \ \ +t_3t_{14}t_{12}t_{234}+q^2t_2t_{14}t_{34}t_{123}+qt_4t_{23}t_{234}-\overline{q}t_3t_{14}t_{134}+q^3t_1t_{23}t_{123}  \nonumber \\
&\ \ \ \ -q^3t_2t_{14}t_{124} +\overline{q}t_{14}t_{12}t_{24}+\overline{q}t_{14}t_{13}t_{34}+(1-q^2)t_2t_4t_{14}t_{12} \nonumber  \\
&\ \ \ \ -\overline{q}^2t_{14}^2+q^2t_{23}^2+(q-\overline{q})t_1t_4t_{14}+q^2(t_1^2+t_4^2-\alpha^2).  \nonumber
\end{align}
\end{prop}

\begin{proof}
The first identity is a repeat of (\ref{eq:13-24}). The third is a result of \cite{BP00}. It can be deduced from (\ref{eq:12-23}), (\ref{eq:top-times-down}); see \cite{Ch22} Example 4.6 for an elegant deduction.

The second identity is obtained by combining (\ref{eq:24-134}) with (\ref{eq:14-234}), (\ref{eq:34-124}) and
\begin{align}
t_{14}t_{34}=\overline{q}t_{1\overline{4}34}+qt_{13}+t_1t_3+t_4t_{134}.   \label{eq:14-34}
\end{align}

To deduce (\ref{eq:123-234}), we first compute
$$t_{123}t_{234}=qt_{1234\overline{3}\overline{2}}+\overline{q}t_{14}+t_1t_4+t_{23}t_0,$$
and then replace $t_{1234\overline{3}\overline{2}}$ via (\ref{eq:12-23-34}).

The identity (\ref{eq:123-134}) is obtained by combing
\begin{align*}
t_{123}t_{134}&=qt_{234\overline{3}}+\overline{q}t_{12\overline{1}4}+t_2t_4+t_{13}t_0,  \\
t_{12\overline{1}4}&=qt_{12}t_{14}-q^2t_{24}-qt_2t_4-qt_1t_{124},  \\
t_{234\overline{3}}&=\overline{q}t_{23}t_{34}-\overline{q}^2t_{24}-\overline{q}t_2t_4-\overline{q}t_3t_{234}.
\end{align*}

\begin{figure}[h]
  \centering
  \includegraphics[width=12cm]{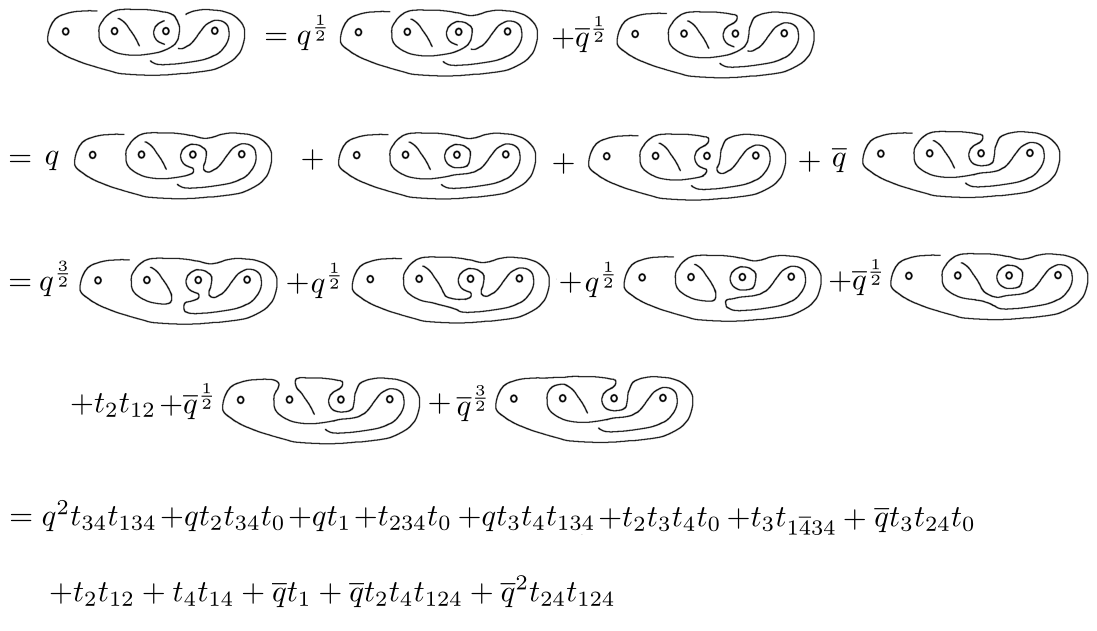}\\
  \caption{Computing $t_{23}t_{12\overline{4}34}$.}\label{fig:computation}
\end{figure}

To deduce (\ref{eq:23-34-124}), we first compute
\begin{align*}
t_{23}t_{12\overline{4}34}=\ &(qt_2t_{34}+\overline{q}t_3t_{24}+t_2t_3t_4+t_{234})t_0+(q^2t_{34}+qt_3t_4)t_{134}   \\
&+(\overline{q}t_2t_4+\overline{q}^2t_{24})t_{124}+\alpha t_1+t_4t_{14}+t_2t_{12}+t_3t_{1\overline{4}34},
\end{align*}
as shown in Figure \ref{fig:computation}, and use (\ref{eq:14-34}) to substitute $t_3t_{1\overline{4}34}$ with elements of $\mathcal{T}_n$.
Then (\ref{eq:23-34-124}) follows from (\ref{eq:34-124}).

Replacing $2$, $3$ in (\ref{eq:top-times-down}) respectively with $23$, $4$ and taking the mirror, we obtain
\begin{align*}
t_{14}t_{1234\overline{3}\overline{2}}&=t_0^2+(\overline{q}t_1t_{234}+qt_4t_{123}+t_1t_4t_{23})t_0+q^2t_{123}^2+\overline{q}^2t_{234}^2  \\
&\ \ \ \ +\overline{q}t_4t_{23}t_{234}+qt_1t_{23}t_{123}+t_{23}^2+t_1^2+t_4^2-\alpha^2.
\end{align*}
Then applying (\ref{eq:12-23-34}) leads to the last identity.
\end{proof}

Use {\it reduction relations} to name these relations and the ones obtained from cyclically permuting indices in
(\ref{eq:24-134'})--(\ref{eq:23-34-124}) as well as the mirror of (\ref{eq:123-234}).


\subsection{Statement and proof}

\begin{thm}
The Kauffman bracket skein algebra of $\Sigma_{0,5}$ 
is generated by $t_1,t_2,t_3,t_4$, $t_0$, $t_{12},t_{13},t_{14},t_{23},t_{24},t_{34}$, $t_{123},t_{124},t_{134},t_{234}$, and the ideal of defining relations is generated by the commuting relations of type $[2,2]$, $[2,3]$, the reduction relations, along with the centralities of $t_1,t_2,t_3,t_4,t_0$.
\end{thm}

\begin{proof}
Given a monomial $\mathfrak{g}$, call it {\it reduced} if it is a product of non-central generators. 
When ${\rm md}_{\mathfrak{g}}(i)=0$ for some $i$, the result on $\mathcal{S}_3$ can be applied.
So we can assume ${\rm md}_{\mathfrak{g}}(i)>0$ for all $i$.

The strategy is to show that using the known relations, any reduced monomial can be converted into a linear combination of certain 
{\it distinguished} ones. The distinguished monomials can be seen to be linearly independent from their leading multicurves; (if it is not so, a new relation would be found).

Up to symmetries, it suffices to consider the cases in the following table.

\begin{center}
\begin{tabular}{|c|c|}
  \hline
  multidegree & distinguished monomials (leading multicurves) \\
  \hline
  $(1,1,1,1)$ & $t_{12}t_{34}$,\ $t_{14}t_{23}$  \\
  \hline
  $(1,1,1,2)$ & $t_{14}t_{234}$ ($t_{2341\overline{4}}$),\ $t_{34}t_{124}$ ($t_{12\overline{4}34}$) \\
  \hline
  $(1,1,2,2)$ & $t_{12}t_{34}^2$,\ $t_{14}t_{23}t_{34}$ \\ 
  \hline
  $(1,2,1,2)$ & $t_{12}t_{24}t_{34}$ ($t_{12\overline{4}34\overline{2}}$), \ $t_{14}t_{23}t_{24}$ ($t_{1\overline{4}23\overline{2}4}$) \\
  \hline
  $(1,1,1,3)$ & $t_{14}t_{24}t_{34}$ \\ 
  \hline
  $(1,2,2,2)$ & $t_{12}t_{34}t_{234}$ ($t_{34}t_{1234\overline{2}}$),\ $t_{14}t_{23}t_{234}$ ($t_{23}t_{1\overline{4}234}$)  \\
  \hline
  $(1,1,2,3)$ & $t_{34}^2t_{124}$ ($t_{12\overline{4}\overline{3}434}$),\ $t_{14}t_{34}t_{234}$ ($t_{1\overline{4}\overline{3}4234}$)  \\
  \hline
  $(1,2,1,3)$ & $t_{24}t_{34}t_{124}$ ($t_{124\overline{2}\overline{4}34}$),\ $t_{14}t_{24}t_{234}$ ($t_{1\overline{4}\overline{2}4234}$)  \\
  \hline
  $(2,2,2,2)$ & $t_{12}^2t_{34}^2$, $t_{14}^2t_{23}^2$  \\
  \hline
  $(2,2,1,3)$ & $t_{14}^2t_{23}t_{24}$ ($t_{14\overline{1}\overline{4}23\overline{2}4}$),\ $t_{12}t_{14}t_{24}t_{34}$ ($t_{124\overline{2}\overline{1}\overline{4}34}$) \\
  \hline
  $(1,1,3,3)$ & $t_{12}t_{34}^3$,\ $t_{14}t_{23}t_{34}^2$ ($t_{23}t_{1\overline{4}\overline{3}234}$) \\
  \hline
  $(1,3,1,3)$ & $t_{12}t_{34}t_{24}^2$ ($t_{124\overline{2}}t_{2\overline{4}34}$),\ $t_{14}t_{23}t_{24}^2$ ($t_{23\overline{2}4}t_{1\overline{4}24}$) \\
  \hline
  $(1,1,2,4)$ & $t_{14}t_{24}t_{34}^2$  \\
  \hline
  $(1,2,1,4)$ & $t_{14}t_{24}^2t_{34}$  \\
  \hline
  $(2,2,2,3)$ & $t_{12}t_{34}^2t_{124}$ ($t_{12}t_{12\overline{4}\overline{3}434}$),\ $t_{14}^2t_{23}t_{234}$ ($t_{23}t_{14\overline{1}\overline{4}234}$)  \\
  \hline
  $(1,2,2,4)$ & $t_{24}t_{34}^2t_{124}$ ($t_{2\overline{4}34}t_{12\overline{4}34}$),\
  $t_{14}t_{24}t_{34}t_{234}$ ($t_{24}t_{1\overline{4}34234}$)   \\
  \hline
  $(1,2,4,2)$ & $t_{23}t_{34}^2t_{123}$ ($t_{234\overline{3}}t_{1234\overline{3}}$),\
  $t_{23}^2t_{34}t_{134}$ ($t_{234\overline{3}}t_{1\overline{3}234}$)   \\
  \hline
  $(1,2,3,3)$ & $t_{14}t_{23}t_{34}t_{234}$  \\
  \hline
  $(1,3,2,3)$ & $t_{12}t_{24}t_{34}t_{234}$ ($t_{12342\overline{4}34\overline{2}}$), \
  $t_{14}t_{23}t_{24}t_{234}$ ($t_{1\overline{4}2\overline{3}\overline{2}4234}$)  \\
  \hline
  $(1,1,3,4)$ & $t_{34}^3t_{124}$ ($t_{12\overline{4}\overline{3}\overline{4}3434}$),\ $t_{14}t_{34}^2t_{234}$ ($t_{34}^2t_{1\overline{4}234}$)  \\
  \hline
  $(1,3,1,4)$ &  $t_{14}t_{24}^2t_{234}$   \\
  \hline
  $(2,2,3,3)$ & $t_{12}^2t_{34}^3$,\ $t_{14}^2t_{23}^2t_{34}$  \\
  \hline
  $(2,3,2,3)$ & $t_{12}^2t_{24}t_{34}^2$ ($t_{12\overline{4}\overline{3}434\overline{2}\overline{1}2}$),\
  $t_{14}^2t_{23}^2t_{24}$ ($t_{14\overline{1}\overline{4}232\overline{3}\overline{2}4}$) \\
  \hline
  $(2,2,2,4)$ & $t_{12}t_{14}t_{24}t_{34}^2$ ($t_{12\overline{4}34}^2$),\ $t_{14}^2t_{23}t_{24}t_{34}$ ($t_{1\overline{4}234}^2$) \\
  \hline
\end{tabular}
\end{center}

As a supplement, if the multidegree is $(1,a,b,c)$ with $a+b+c=9$, then each relation is implied by the above ones. To see this, note that if $S$ is a simple curve with ${\rm md}_S=(1,a,b,c)$ with $a+b+c=9$, then any two of admissible expressions are congruent through a shortenable arc.

Now that all possible monomials have been checked, the proof is complete.
\end{proof}

\bigskip

\noindent
Haimiao Chen (orcid: 0000-0001-8194-1264)\ \ \  \emph{chenhm@math.pku.edu.cn} \\
Department of Mathematics, Beijing Technology and Business University, \\
Liangxiang Higher Education Park, Fangshan District, Beijing, China.

\end{document}